\documentclass{article}
\usepackage[T2A]{fontenc}
\usepackage[cp1251]{inputenc}
\usepackage[english,russian]{babel}
\usepackage[tbtags]{amsmath}
\usepackage{amsfonts,amssymb,mathrsfs,amscd}

\usepackage[hyper]{msb-a}

\makeatletter
\gdef\No{{\select@language{russian}\textnumero}}
\makeatother

\JournalName{}
\numberwithin{equation}{section}
\theoremstyle{plain}
\newtheorem{theorem}{Теорема}
\newtheorem{lemma}{Лемма}[section]
\newtheorem{propos}{Предложение}
\theoremstyle{definition}
\newtheorem{definition}{Определение}
\newtheorem{proof}{Доказательство}
\newtheorem{remark}{Замечание}

\def\CC{\mathbb C}

\def\RS{\mathfrak R}
\def\bz{\mathbf z}
\def\sH{\mathscr H}

\graphicspath{{figures/}}

\def\R{{\mathbb R}}
\def\N{{\mathbb N}}
\def\Z{{\mathbb Z}}
\def\dd{{\rm d}}

\newcommand{\SEtwo}{\operatorname{SE}(2)}
\newcommand{\PTRtwo}{\operatorname{PT}(\R^2)}
\newcommand{\sign}{\operatorname{sign}}
\newcommand{\am}{\operatorname{am}}
\newcommand{\cn}{\operatorname{cn}}
\newcommand{\sn}{\operatorname{sn}}
\newcommand{\dn}{\operatorname{dn}}
\newcommand{\EE}{\operatorname{E}}
\newcommand{\sech}{\operatorname{sech}}
\renewcommand{\dd}{\operatorname{d}\nolimits}

\newcommand{\Exp}{\operatorname{Exp}\nolimits}
\newcommand{\M}{\mathbb{M}}
\newcommand{\s}{\mathrm{s}}

\newcommand{\tx}{\tilde{x}}
\newcommand{\ty}{\tilde{y}}
\newcommand{\tth}{\tilde{\theta}}
\renewcommand{\i}{ \boldsymbol{i}}
\newcommand{\D}{U}
\begin{document}

\title{Задача быстродействия на группе движений плоскости  с управлением в полукруге}
\author[A.\,P.~Mashtakov]{А.\,П.~Маштаков}
\address{Институт программных систем им.~А.\,К.~Айламазяна РАН}
\email{alexey.mashtakov@gmail.com}

\date{\today}
\udk{517.977}

\maketitle

\begin{fulltext}

\begin{abstract}
Исследуется задача быстродействия на группе движений плоскости с управлением в полукруге. Рассматриваемая управляемая система задает модель машины на плоскости, которая может двигаться вперед и вращаться на месте. Оптимальные по заданной внешней стоимости траектории такой системы используются в обработке изображений 
для поиска выделяющихся кривых. В частности, такие траектории используются в анализе медицинских изображений при поиске сосудов на фото сетчатки глаза человека. Задача представляет интерес в геометрической теории управления, как модельный пример, в котором множество значений управляющих параметров содержит ноль на границе. В работе изучен вопрос управляемости и существования оптимальных траекторий. 
На основе анализа 
гамильтоновой системы принципа максимума Понтрягина 
найден явный вид экстремальных управлений и траекторий. Частично исследован вопрос оптимальности экстремалей. 
Описана структура оптимального синтеза. 

Библиография: 33 названия.
\end{abstract}


\markright{Задача быстродействия на $\SEtwo$  с управлением в полукруге}

\footnotetext[0]{Исследование выполнено при финансовой поддержке РНФ в рамках научного проекта \No~17-11-01387-П в Институте программных систем им. А.К. Айламазяна Российской академии наук}

\section{Введение}

В 1957 году Л.~Дубинс описал в работе~\cite{dubinscurves1957} задачу поиска кратчайших путей для автомобиля (машины), движущегося по плоскости из начальной конфигурации (положение и направление) в конечную конфигурацию. В постановке Дубинса автомобиль может двигаться только вперед (не имеет задней передачи), а его угловая скорость ограничена. Позже, в 1990 году, Дж.~Ридс и Л.~Шепп рассмотрели в работе~\cite{reedsoptimal1990} ту же задачу, но применительно к автомобилю, у которого есть возможность движения назад. В обеих статьях основное внимание уделяется описанию и доказательству общей формы кратчайших путей без предоставления явных решений для заданных граничных условий (начальной и конечной конфигураций). 

Важным недостатком моделей Дубинса и Ридса-Шеппа является то, что траектории движения автомобиля имеют ограниченную кривизну. Это означает, что машина не может вращаться на месте. В то же время, в робототехнике принято рассматривать модель автомобиля с независимыми приводами для двух колес. В такой системе два колеса могут вращаться в разных направлениях с одинаковой скоростью, что обеспечивает вращение на месте. Таким образом, естественное обобщение машины Ридса-Шеппа приводит к модели автомобиля, траектории которого имеют неограниченную кривизну. Оптимальный синтез в такой системе получен Ю.Л.~Сачковым~\cite{yuriSE2}. Аналогичное обобщение машины Дубинса приводит к модели, предложенной Р. Дайтсом~\cite{Duits18}.

Рассмотрим модель автомобиля, движущегося по плоскости, см. Рис.~\ref{fig:Intro}. Автомобиль имеет два параллельных колеса, равноудаленных от центра масс, который совпадает с серединой оси колесной пары. Оба колеса имеют независимые приводы, которые могут вращаться вперед и назад, так что соответствующее качение колес происходит без проскальзывания. Конфигурация системы описывается тройкой $ q = (x, y, \theta) \in \M = \R^ 2 \times S^1 $, где $ (x, y) \in \R^2 $ --- центральная точка, а $\theta \in S^1 $ --- угол ориентации автомобиля, совпадающий с направлением колес. Заметим, что конфигурационное пространство $\M$ образует группу Ли $\SEtwo$ --- группу движений плоскости:
$$\M = \R^2 \times S^1\simeq\SEtwo.$$

\begin{figure}[h]
\centering
\begin{minipage}{0.68\hsize}
\includegraphics[width=0.9\hsize]{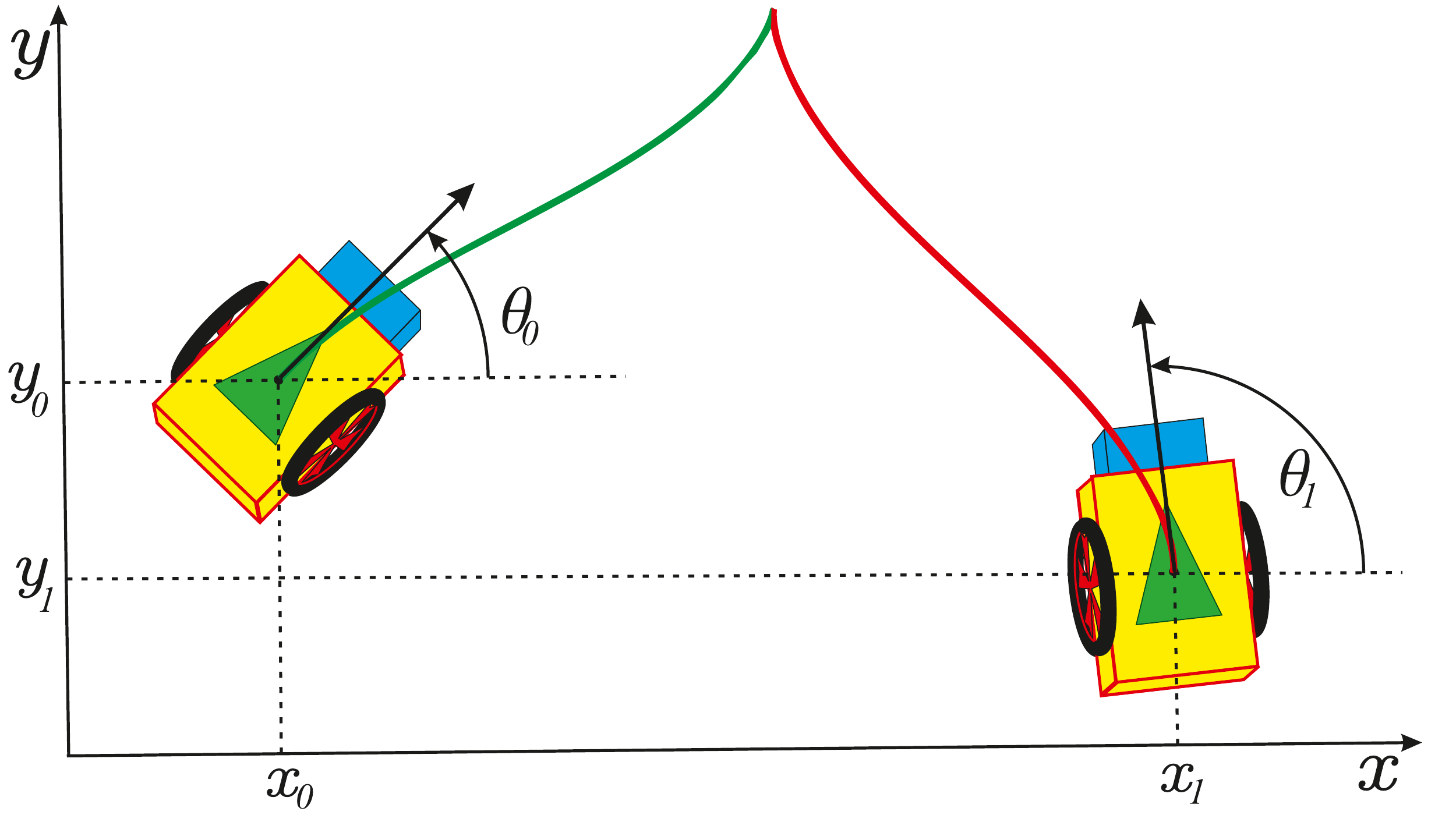}
\end{minipage}
\begin{minipage}{0.31\hsize}
\includegraphics[width=0.99\hsize]{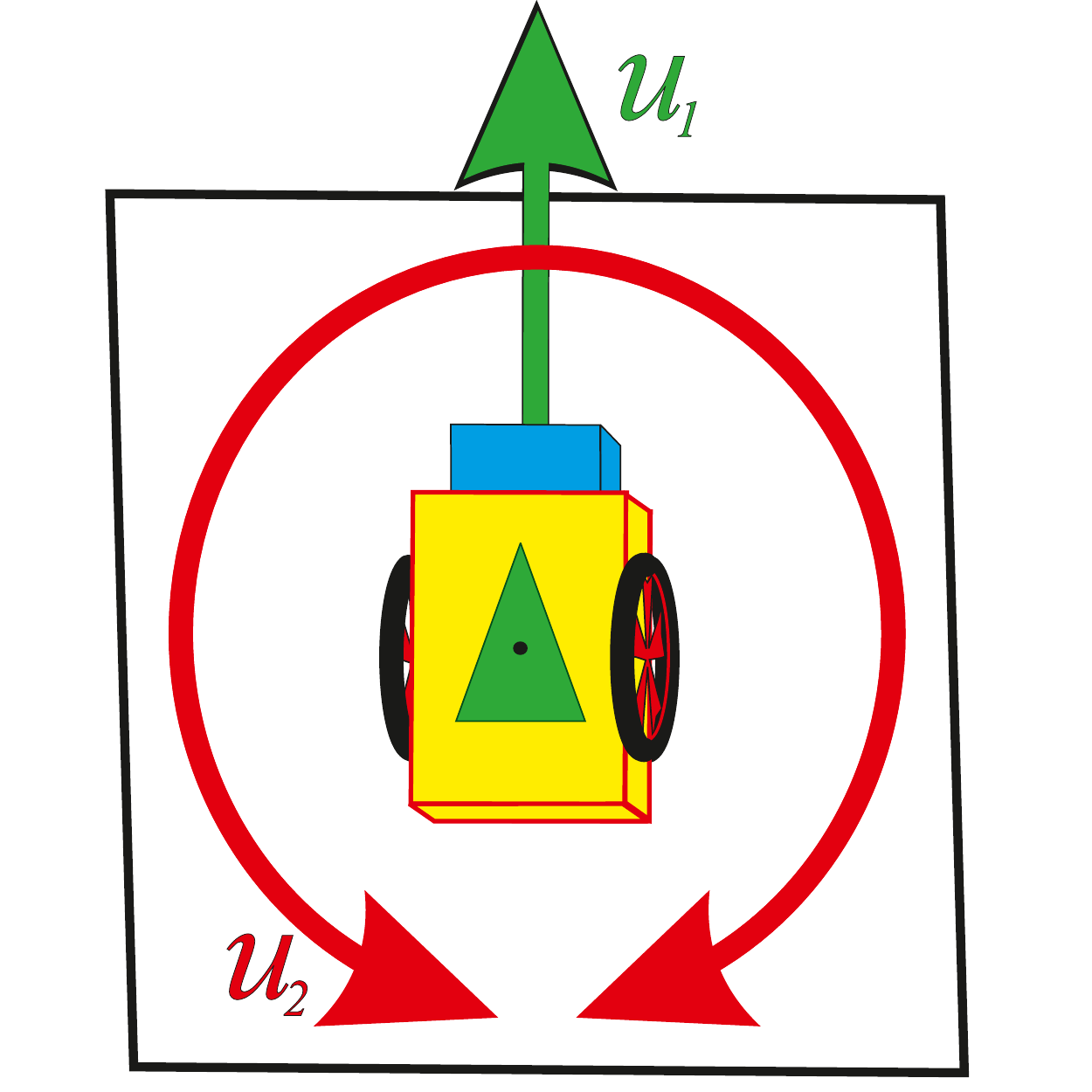}
\end{minipage}
\caption{\textbf{Слева:} 
Конфигурация системы ``машина на плоскости'' определяется положением $(x,y) \in \R^2$ и углом ориентации $\theta \in S^1$.
В модели Сачкова машина может двигаться вперед (часть траектории выделена зеленым) и назад (выделена красным). В момент смены направления движения траектория имеет точку возврата. \textbf{Справа:} 
Машина имеет два управления: $u_1$, отвечающее за движение вперед, и $u_2$, отвечающее за повороты. В модели Дайтса движение назад запрещено~$u_1 \geq 0$.
}
\label{fig:Intro}
\end{figure}

С точки зрения водителя у автомобиля есть два управления: акселератор $u_1$ и рулевое колесо $u_2$. В таких обозначениях динамика автомобиля описывается следующей управляемой системой, см.~\cite{Laumond1986}:
\begin{equation}\label{eq:contsys}
\begin{cases}
\dot{x} = u_1 \cos \theta,\\
\dot{y} = u_1 \sin \theta, \\
\dot{\theta} = u_2,
\end{cases}
\qquad 
\begin{array}{l} 
(u_1,u_2) \in U \subset \R^2,\\
(x,y,\theta) \in \M = \R^2 \times S^1.
\end{array}
\end{equation}

Задача заключается в поиске управления $(u_1(t), u_2(t))$ и соответствующей траектории $\gamma (t) = (x(t), y(t), \theta(t)) $
на временном интервале $ t \in [0, T] $, удовлетворяющей системе~(\ref{eq:contsys}) с граничными условиями
\begin{equation}\label{eq:boundsint}
\gamma(0) = (x_0, y_0, \theta_0) = q_0, \qquad \gamma(T) = (x_1, y_1, \theta_1) = q_1, 
\end{equation}
где $q_0, q_1 \in \M $ --- это начальная и конечная конфигурации системы.

Такая задача в общем случае допускает бесконечное число решений. Выбор конкретного решения традиционно осуществляется путем введения критерия качества различных траекторий. Для этого к системе добавляется функционал стоимости, который требуется минимизировать, чтобы найти оптимальное решение. Классической 
 является задача быстродействия 
\begin{equation}\label{eq:intL}
T \to \min. 
\end{equation}
Еще один естественный критерий оптимальности --- минимизация ``маневра''
\begin{equation}\label{eq:intLSR}
\int_0^T \sqrt{u_1^2 (t) + u_2^2(t)} \dd t \to \min, 
\end{equation}
где конечное время $T \geq 0 $ свободно. 

Эти два критерия могут быть эквивалентны при выборе подходящей области допустимых управлений. Например, в случае $u_1^2 (t) + u_2^2 (t) \leq 1$ два критерия эквивалентны, поскольку любая нетривиальная оптимальная траектория может быть натурально параметризована $u_1^2 (t) + u_2^2 (t) = 1$.
 
\begin{remark}
Задача управления~(\ref{eq:contsys}),~(\ref{eq:boundsint}),~(\ref{eq:intLSR}) при $U=-U$ является задачей оптимизации кривой в пространстве $\M$, снабженном естественной метрикой, определенной для кривых $\gamma(t) = (x(t), y(t), \theta(t))$, удовлетворяющих условию пропорциональности вектора $ (\dot{x}(t), \dot{y}(t)) $ вектору $ (\cos \theta (t) , \sin \theta (t)) $. Формулируя задачу таким образом, она становится одним из простейших примеров субримановой (СР) геометрии~\cite{montgomery}: касательный вектор $ \dot \gamma(t) $ должен лежать в плоскости, натянутой на вектора $ (\cos \theta (t), \sin\theta (t), 0) $ и $ (0,0,1)$.
\end{remark}

Различные множества допустимых управлений $U \ni (u_1, u_2)$ приводят к разным моделям автомобиля на плоскости, см.~\cite{sussman}. Например, см. Рис.~\ref{fig:U}, задача минимизации времени для автомобиля с управлением
\begin{itemize}
\item $u_1 = 1$, $|u_2| \leq \kappa$, $ \kappa> 0 $ приводит к модели Дубинса~\cite{dubinscurves1957};
\item $|u_1| = 1 $, $|u_2| \leq \kappa $, $ \kappa> 0 $ приводит к модели Ридса-Шеппа~\cite{reedsoptimal1990};
\item $u_1^2+u_2^2 \leq 1$ приводит к модели машины, траектории которой задаются субримановыми кратчайшими, исследованными Сачковым~\cite{yuriSE2};
\item $u_1^2+u_2^2 \leq 1$, $u_1 \neq 0$ приводит к модели, исследованной Берестовским~\cite{Berestovskiy};
\item $u_1^2+u_2^2 \leq 1$, $u_1 \geq 0$ приводит к модели машины, движущейся вперед и поворачивающейся на месте, предложенной Дайтсом~\cite{Duits18}.
\end{itemize}

\begin{figure}[h]
\begin{minipage}{0.99\hsize}
\includegraphics[width=0.9\hsize]{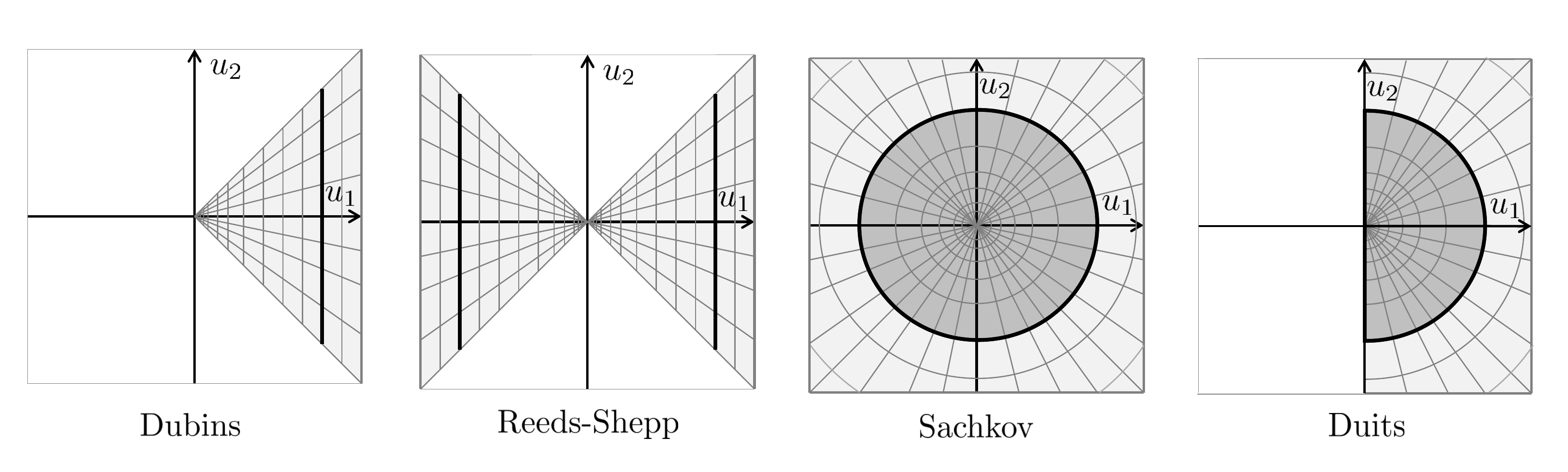}
\end{minipage}
\caption{Множество допустимых управлений для различных моделей.
}
\label{fig:U}
\end{figure}

В данной работе исследуется модель Р. Дайтса~\cite{Duits18}. Актуальность модели обусловлена приложением к обработке изображений. Оптимальные для заданной внешней стоимости траектории такой системы используются для поиска выделяющихся кривых. В частности, такие траектории используются в анализе медицинских изображений при поиске сосудов на фото сетчатки глаза. Модель Дайтса была разработана, чтобы устранить проблему точек возврата, возникающую в методе трассировки сосудов с помощью субримановых (СР) кратчайших~\cite{CIARP}. Отметим, что для устранения точек возврата применялся также другой подход~\cite{PTR2}, основанный на замене конфигурационного пространства с $\SEtwo$ на $\PTRtwo$. Такой подход позволил существенно сократить число ситуаций, в которых возникают точки возврата, однако полностью проблема устранена не была. Более подробно применение модели в обработке изображений обсуждается в следующем разделе.

Помимо актуальности с точки зрения приложений рассматриваемая задача представляет самостоятельный интерес в геометрической теории управления~\cite{notes}, как модельный пример задачи оптимального управления, в которой нулевое управление находится на границе множества управляющих параметров. Наиболее общий подход~\cite{LokTrig} к решению родственных задач основан на выпуклой тригонометрии. Подход покрывает класс задач оптимального управления с двумерным управлением, принадлежащим произвольному выпуклому компакту, содержащему ноль внутри. Обобщить этот подход на случай, когда ноль лежит на границе, непосредственно не удается. Для систем такого типа требуется развитие новых методов. В данной статье детально исследован частный случай такой системы.


Работа имеет следующую структуру. Помимо введения, в котором описана история задачи и отражена ее актуальность в мобильной робототехнике и теории управления, статья содержит пять разделов и заключение. Раздел~\ref{sec:Imaging} посвящен приложениям в моделировании биологических зрительных систем и обработке изображений. В разделе~\ref{sec:Statement} приводится формальная постановка задачи. В разделе~\ref{sec:Exist} доказывается полная управляемость и существование оптимальных управлений. В разделе~\ref{sec:PMP} к задаче применяется принцип максимума Понтрягина (ПМП), 
 исследуется гамильтонова система ПМП 
и приводится явный вид экстремальных управлений и траекторий. В разделе~\ref{sec:optimal} исследуется оптимальность экстремальных траекторий и структура оптимального синтеза.
В заключении подводится итог работы.

\section{Приложение в моделях зрения и обработке изображений}\label{sec:Imaging}

Принципы устройства биологических зрительных систем вызывают большой интерес среди исследователей во многих областях науки. Важным направлением является разработка и изучение реалистичных математических моделей, описывающих определенный этап обработки зрительного сигнала. Математическая модель первичной зрительной коры мозга как субриманова структура в пространстве позиций и направлений была предложена Ж.~Петито~\cite{Petitot_B}. Затем модель была уточнена Дж.~Читти и А.~Сарти~\cite{Citti_Sarti_B}. 

Модель Петито-Читти-Сарти представляет первичную зрительную кору как субриманову структуру на группе Ли $\SEtwo$. 
Иллюзорные контуры в такой модели определяются как траектории системы~(\ref{eq:contsys}), на которых достигается минимум функционала~(\ref{eq:intL}) на множестве управлений $u_1^2+u_2^2 \leq 1$. 
Согласно этой модели процесс дополнения скрытых контуров происходит путем минимизации энергии возбуждения нейронов, воспринимающих визуальную информацию в областях повреждения. Такой процесс моделируется действием оператора гипоэллиптической диффузии, изученного в~\cite{Gauthier,Gauthier2}. Результирующие кривые являются субримановыми кратчайшими. Такие кривые применяются для восстановления поврежденных изображений~\cite{OptimalInpainting} и используются для объяснения некоторых зрительных иллюзий~\cite{MashtakovDGA}. 


В работах~\cite{DuitsJMIV2014,Boscain2014} указано, что не все траектории модели Петито-Читти-Сарти согласованы с экспериментами из психологии зрения по построению поля ассоциаций, описанными в работе~\cite{AssocField}. Поле ассоциаций представляет собой множество граничных условий (позиций и направлений), соединяемых иллюзорным контуром в процессе работы зрительной системы человека. Экспериментально было установлено, что не любые граничные условия соединяются. Для тех граничных условий, которые соединяются, иллюзорный контур является гладкой кривой на плоскости $\R^2$. Заметим, что в модели Петито-Читти-Сарти любые заданные граничные условия могут быть соединены оптимальной траекторией. При этом проекция оптимальной траектории на плоскость может иметь точки возврата (не является гладкой). 

Естественным уточнением модели Петито-Читти-Сарти является сужение задачи на класс граничных условий, соединяемых оптимальной траекторией без точек возврата. Исследование множества таких граничных условий произведено в~\cite{DuitsJMIV2014,Boscain2014}. Перспективным направлением является модификация модели таким образом, чтобы траектории с точками возврата были априори исключены. Такая модификация дается моделью Дайтса~\cite{Duits18}.  

 Принципы работы биологических зрительных систем активно используются в компьютерном зрении. На основе этих принципов создаются эффективные методы обработки изображений: улучшение, сегментация, восстановление, поиск объектов. Так, например, в работе~\cite{Duits2007} предложен подход к обработке изображений, основанный на поднятии изображения в расширенное пространство позиций и направлений (invertible orientation scores). После такого подъема, субримановы (СР) кратчайшие используются для поиска выделяющихся кривых~\cite{SIAM, SO3JMIV, SE3JDCS}. 

\begin{figure}[h]
\centering
\includegraphics[width=0.387\hsize]{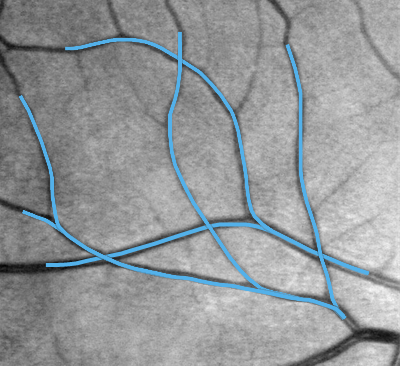}$\quad$
\includegraphics[width=0.39\hsize]{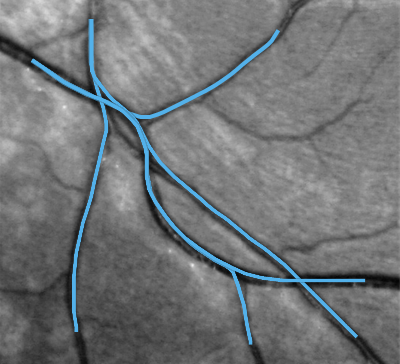}
\caption{Трассировка сосудов на фото сетчатки глаза с помощью траекторий на $\SEtwo$ с управлением в полукруге. Адаптировано из~\cite{Scharpach}.}
\label{fig:Vessel}
\end{figure}


В частности, задача поиска выделяющихся кривых возникает в анализе медицинских изображений при поиске кровеносных сосудов на фото сетчатки глаза человека, см. Рис.~\ref{fig:Vessel}. Типичные проблемы для промышленных систем трассировки возникают в ситуациях, когда сосуды на фото пересекаются. Решение этой проблемы достигается путем поднятия изображения в расширенное пространство $\M$ позиций и направлений. Метод трассировки сосудов с помощью СР кратчайших на $\M$ предложен в~\cite{CIARP}.

 Недостатком СР модели является наличие точек возврата. 
Такие кривые нежелательны в задаче трассировки сосудов. Ограничение управления на полукруг устраняет этот недостаток. При этом точки возврата становятся точками поворота на месте. Возникновение точки поворота типично наблюдается в местах ветвления сосудов. Более подробно с результатами трассировки сосудов с помощью оптимальных траекторий на $\SEtwo$ с управлением в полукруге можно ознакомиться в~\cite{Scharpach}.      

\section{Постановка задачи}\label{sec:Statement}
Движением метрического пространства называется преобразование, сохраняющее расстояние между точками.  В данной работе рассматривается группа $\SEtwo$ собственных евклидовых движений плоскости $\R^2$. Любое собственное движение $\R^2$ представляется композицией параллельного переноса и поворота плоскости. Таким образом, элемент $\SEtwo$ задается тройкой чисел $(x,y,\theta)$,
где $(x, y)\in\R^2$ --- вектор параллельного переноса, $\theta\in S^1=\R / 2 \pi \Z$  --- угол поворота.

Рассматривается следующая управляемая система:
\begin{equation}\label{eq:contsyst}
\left\{\begin{array}{l}
\dot{x} = u_1 \cos{\theta}, \\[0.005\linewidth]
\dot{y} = u_1 \sin{\theta}, \\[0.005\linewidth]
\dot{\theta} = u_2, 
\end{array}\right.
\quad
\begin{array}{l}
(x,y,\theta) =q \in \SEtwo=\M, \\[0.01\linewidth]
u_1^2+u_2^2\leq1,\, u_1\geq0.
\end{array}
\end{equation}
Исследуется задача быстродействия: по заданным граничным условиям $q_0$, $q_1 \in \M$ требуется найти управления $u_1(t)$, $u_2(t)$ такие, что соответствующая траектория $\gamma:[0,T] \to \M$ переводит систему из начального состояния $q_0$ в конечное состояние $q_1$ за минимальное время
\begin{equation}\label{eq:bounds}
\gamma(0)= q_0,\quad \gamma(T)=q_1, \qquad T\to\min.
\end{equation}
В данной постановке управления $u_i$ принадлежат классу $L^\infty([0,T],\R)$, а соответствующие траектории $\gamma$ являются  липшицевыми кривыми на $\M$.
\begin{remark}\label{rem:bounds}
Система~(\ref{eq:contsyst}) инвариантна относительно параллельных переносов и поворотов плоскости. В силу этого без ограничения общности можно свести исследование при произвольном $q_0=(x_0,y_0,\theta_0)$ к случаю $q_0=(0,0,0).$
\end{remark}
\begin{remark}\label{rem:u}
Управление  $u_1$ задает инфинитезимальный параллельный перенос в направлении $\theta$, а $u_2$ --- поворот плоскости $\R^2_{x,y}$. При $u_2>0$ поворот происходит от оси $O_x$ к $O_y$, а при $u_2<0$ --- в обратном направлении.   
\end{remark}

Классический подход~\cite{notes} к задаче~(\ref{eq:contsyst})--(\ref{eq:bounds}) состоит из следующих этапов:
\begin{enumerate}
\item \label{enum:st1} доказательство существования решения;
\item \label{enum:st2} параметризация экстремалей $\gamma$ с помощью принципа максимума Понтрягина (ПМП);
\item \label{enum:st3}
определение времени разреза $t_{cut}(\gamma)$, после которого $\gamma$ теряет оптимальность;
\item \label{enum:st4} выбор оптимальной траектории для заданных граничных условий.
\end{enumerate}

\begin{remark}\label{rem:uniq}
Заметим, что при $t\in[0,t_{cut})$ оптимальная траектория~$\gamma$ является единственной. В конечную точку $\gamma(t_{cut})$ могут приходить нес\-колько оптимальных траекторий, в том числе бесконечно много. 
\end{remark}

Далее в разделах~\ref{sec:Exist}---\ref{sec:PMP} полностью проводятся этапы~\ref{enum:st1} и~\ref{enum:st2} соответственно. В разделе~\ref{sec:optimal} обсуждаются этапы~\ref{enum:st3} и~\ref{enum:st4}, которые на данный момент являются открытыми математическими задачами в общем случае. Проведение этих этапов позволит построить оптимальный синтез в задаче. На данный момент неизвестно время разреза $t_{cut}$  для всех возможных экстремалей $\gamma$. Сложность его вычисления заключается в решении трансцендентных уравнений, возникающих при параметризации экстремалей. 


\section{Существование решения}\label{sec:Exist}
При исследовании задачи~(\ref{eq:contsyst})--(\ref{eq:bounds}) возникает вопрос существования допустимой траектории, соединяющей граничные условия~(\ref{eq:bounds}). Если при любых $q_0$, $q_1 \in \M$ ответ положителен, то управляемая система называется вполне управляемой. Докажем конструктивно полную управляемость системы~(\ref{eq:contsyst}).

В силу Замечания~\ref{rem:bounds} положим $q_0=(0,0,0)$. Обозначим $q_1=(x_1,y_1,\theta_1)$. Искомое управление будем строить из трех частей: вращение на месте на угол $\alpha = \arg(x_1+\i y_1)\in (-\pi,\pi]$, движение вперед на величину $l=\sqrt{x_1^2+y_1^2}$, вращение на месте на угол $\beta = (\theta_1-\alpha)\mod 2\pi \in (-\pi,\pi]$. Здесь $\arg(a + \i b)$ --- аргумент комплексного числа $a+\i b$.  
Введем следующие обозначения: 
$$s_1 = \sign\alpha, \;s_2 = \sign\beta; \quad
t_0^1 = |\alpha| \in [0, \pi], \;t_0^2 = |\alpha| +l, \;T = |\alpha| +l + |\beta|.$$

Искомое управление $\bar{\mathbf{u}}=(\bar{u}_1,\bar{u}_2)$ имеет вид 
\begin{equation*}
\bar{\mathbf{u}}(t)=\left\{
\begin{array}{l l}
(0,s_1), & \text{при } t\in[0,t_0^1|),\\
(1,0), & \text{при } t\in[t_0^1,t_0^2),\\
(0,s_2), & \text{при } t\in[t_0^2,T].
\end{array}
\right.
\end{equation*}

Соответствующая траектория $\bar{\gamma}$ задается следующей таблицей:
$$\textrm{
\begin{tabular}{ c | c c c }
$t\in         $&$ [0,t_0^1) $&$ [t_0^1, t_0^2)$&$ [t_0^2,T]$ \\
\hline
$\bar{x}(t)         $&$ 0             $&$ (t- s_1 t_0^1) \cos t_0^1  $&$ x_1$ \\
$\bar{y}(t)         $&$ 0             $&$(s_1 t- t_0^1) \sin t_0^1 $&$y_1$ \\
$\bar{\theta}(t)  $&$ s_1 t        $&$ s_1 t_0^1      $&$ s_1 t_0^1 + s_2 (t-t_0^2)$,
\end{tabular}
}$$  

Вычисление $\theta(T)=s_1 t_0^1 + s_2 (t-t_0^2)=\alpha+\beta=\theta_1$ показывает, что траектория $\bar{\gamma}$ удовлетворяет граничным условиям $\bar{\gamma}(0)=(0,0,0)$, $\bar{\gamma}(T)=(x_1,y_1,\theta_1)$. Итак мы построили управление $\bar{\mathbf{u}}$, переводящее систему из $q_0$ в $q_1$, и тем самым доказали свойство полной управляемости системы~(\ref{eq:contsyst}). 

Далее возникает вопрос существования оптимальных траекторий: всегда ли существует допустимая траектория, удовлетворяющая условиям~(\ref{eq:bounds}), на которой достигается минимальное значение $T$? Положительный ответ на этот вопрос дается теоремой Филиппова~\cite{notes,zelikin}.

Подводя итог, в данном разделе доказана следующая теорема:
\begin{theorem}\label{thm:exist}
Решение задачи быстродействия~(\ref{eq:contsyst})--(\ref{eq:bounds}) существует для любых граничных условий $q_0$, $q_1 \in \M$.
\end{theorem}

\section{Принцип максимума Понтрягина}\label{sec:PMP}
Исследуется задача быстродействия 
\begin{equation}\label{eq:contsystu}
\left\{\begin{array}{l l l}
\dot{x} = u_1 \cos{\theta}, & \;\;\; x(0)=0, & x(T)=x_1, \\[0.005\linewidth]
\dot{y} = u_1 \sin{\theta},  & \;\;\; y(0)=0, & y(T)=y_1, \\[0.005\linewidth]
\dot{\theta} = u_2, & \;\;\; \theta(0)=0, & \theta(T)=\theta_1,
\end{array}\right.
\quad
\begin{array}{l}
u_1\geq 0, \\[0.005\linewidth]
u_1^2+u_2^2 \leq 1,
\end{array}\quad
 T \to \min.
\end{equation}

Обозначим через $\D$ множество допустимых управлений --- полукруг $$\D=\{(u_1,u_2)\;|\;u_1\geq 0,\;u_1^2+u_2^2 \leq 1\}.$$

Применим к~(\ref{eq:contsystu}) необходимое условие оптимальности --- ПМП~\cite{notes,zelikin}. 

Введем укороченную функцию Понтрягина
\begin{equation}\label{eq:pontHam}
H_{u} =u_1 (p_1 \cos{\theta} +p_2 \sin{\theta}) + u_2 p_3, \quad (p_1, p_2, p_3) \in T_e^*\M \simeq \R^3.
\end{equation}
Пусть $(u(t), q(t))$, $t \in[0,T]$, --- оптимальный процесс, тогда выполнены условия:
\begin{enumerate}
\item гамильтонова система $\displaystyle \dot{p}=-\frac{\partial H_u}{\partial q},\,\,  \dot{q}=\frac{\partial H_u}{\partial p}$;
\item условие максимума $H_{u(t)} (p(t),q(t)) ={\underset{u \in \D}{\max}}\, H_{u} (p(t),q(t))=H \in \{0,1\}$.
\end{enumerate}
Случай $H=0$ называется анормальным, случай $H=1$ --- нормальным.

Естественными координатами для левоинвариантных задач на группах Ли являются левоинвариантные гамильтонианы, линейные на слоях кокасательного расслоения~\cite{notes}. Обозначим их через $h_1$, $h_2$ и $h_3$:
\begin{equation*}
h_1=p_1 \cos\theta+p_2 \sin\theta, \quad h_2=p_3, \quad h_3 =p_1 \sin\theta-p_2 \cos\theta . 
\end{equation*}

Функция Понтрягина примет вид
\begin{equation*}\label{eq:pontHam1}
H_{u} =u_1 h_1+u_2 h_2.
\end{equation*}

\textbf{Гамильтонова система ПМП} имеет вид:
\begin{equation}\label{eq:hamsys}
\left\{\begin{array}{l}
\dot{x} = u_1 \cos{\theta}, \\[0.005\linewidth]
\dot{y} =u_1 \sin{\theta}, \\[0.005\linewidth]
\dot{\theta} =u_2,
\end{array}\right.
\quad
\left\{\begin{array}{l}
\dot{h}_1 = -u_2 h_3, \\[0.005\linewidth]
\dot{h}_2 = u_1 h_3, \\[0.005\linewidth]
\dot{h}_3 =u_2 h_1.
\end{array}\right.
\end{equation}

Подсистема на переменные состояния $x$, $y$, $\theta$ называется \emph{горизонтальной} частью, а подсистема на сопряженные переменные $h_1$, $h_2$, $h_3$ --- \emph{вертикальной} частью гамильтоновой системы ПМП. Решение вертикальной части вместе с условием максимума определяет экстремальные управления, а решение горизонтальной части --- экстремальные траектории. 

\textbf{Условие максимума} имеет вид $H=\max\limits_{u\in \D} H_u$.

При $h_1<0$ гамильтониан равен $H=|h_2|$, максимум достигается при
\begin{equation}\label{eq:abnormh1min}
u_1=0, \quad u_2=\begin{cases}\sign h_2, \text{ при } h_2\neq 0,
\\
\forall u_2 \in I=[-1,1], \text{ при } h_2=0.\end{cases}
\end{equation}

При $h_1=0$ гамильтониан равен $H=|h_2|$, максимум достигается при
\begin{equation}\label{eq:abnormh1eq}
\begin{cases}u_1=0, \; u_2=\sign h_2, \text{ при } h_2\neq 0,
\\
\forall (u_1, u_2) \in \D, \text{ при } h_2=0.\end{cases}
\end{equation}

При $h_1>0$ гамильтониан равен $H=\sqrt{h_1^2+h_2^2}$, максимум достигается при
\begin{equation}\label{eq:abnormh1max}
 \mathbf{u}=\frac{\mathbf{h}}{\sqrt{h_1^2+h_2^2}} , \quad \mathbf{u}=\left(\begin{array}{c}u_1\\u_2\end{array}\right),\; \mathbf{h}=\left(\begin{array}{c}h_1\\h_2\end{array}\right).
\end{equation}

\begin{remark}\label{rem:hmax}
Выражение экстремального управления при $h_1>0$ имеет следующую геометрическую интерпретацию. Функция $H_u = u_1 h_1 + u_2 h_2=
\left\langle \mathbf{u}, \mathbf{h}\right\rangle$ достигает максимума по $(u_1,u_2)\in \D$, когда векторы $\mathbf{u}$ и $\mathbf{h}$ сонаправлены и евклидова длина $\|\mathbf{u}\|$ имеет максимальное значение, то есть $u_1^2+u_2^2=1$. Таким образом $\mathbf{u}=\mathbf{h}/H$, где $H=\sqrt{h_1^2+h_2^2}$.
\end{remark}

Заключаем, что гамильтониан ПМП имеет вид 
\begin{equation}\label{eq:HamMax}
H = {\underset{u\in \D}{\max}}H_{u}=
\left\{\begin{array}{l l} 
|h_2|,  & \text{при } h_1 \leq 0,\\[0.005\linewidth]
\sqrt{h_1^2+h_2^2},  & \text{при } h_1 > 0.
\end{array} \right.
\end{equation}

Экстремальное управление при $H=1$ задается следующей таблицей:
\begin{equation}\label{eq:uextr}
\textrm{
\begin{tabular}{ c | l l l}
 &$ h_1 < 0$&$ h_1=0$&$ h_1>0$ \\
\hline
$ h_2 < 0  $&$ u_1=0, \,u_2=-1$&$ (u_1, u_2)=(0,-1)$&$ (u_1, u_2)=(h_1,h_2)$ \\
$ h_2 = 0  $&$ u_1=0, \,u_2 \in I$ &$(u_1, u_2)\in \D $ &$(u_1, u_2)=(h_1,h_2)$ \\
$ h_2 > 0  $&$ u_1=0, \,u_2=1 $&$   (u_1, u_2)=(0,1) $&$ (u_1, u_2)=(h_1,h_2)$.
\end{tabular}
}  
\end{equation}

\subsection{Анормальный случай $H=0$}
При $h_1<0$ из~(\ref{eq:HamMax}) следует $h_2=0$. В силу~(\ref{eq:abnormh1min}) имеем $u_1=0$, $u_2 \in I=[-1,1]$. Гамильтонова система
\begin{equation}\label{eq:hamsysh1na}
\left\{\begin{array}{l l}
\dot{x} = 0,  & x(0)=0,\\[0.005\linewidth]
\dot{y} =0,  & y(0)=0,\\[0.005\linewidth]
\dot{\theta} =u_2, & \theta(0)=0,
\end{array}\right.
\quad
\left\{\begin{array}{l l}
\dot{h}_1 = -u_2 h_3,  & h_1(0)=h_{10},\\[0.005\linewidth]
\dot{h}_3 =u_2 h_1,& h_3(0)=h_{30}
\end{array}\right.
\end{equation}
имеет для любого допустимого управления $u_2(t)$ единственное решение
\begin{equation}\label{eq:hamsysh1natr}
\left\{\begin{array}{l }
x(t) = 0,  \\[0.005\linewidth]
y(t) =0,  \\[0.005\linewidth]
\theta(t) =U_2(t),
\end{array}\right.
\quad
\left\{\begin{array}{l}
h_1(t) =  h_{10} \cos U_2(t) - h_{30} \sin U_2(t),\\[0.005\linewidth]
h_3(t) =  h_{30} \cos U_2(t)+h_{10} \sin U_2(t),
\end{array}\right.
\end{equation}
где $$\displaystyle U_2(t)=\int_0^t u_2(\tau) d \tau.$$

При $h_1=0$ из~(\ref{eq:HamMax}) следует $h_2=0$. Поскольку гамильтониан $H=|h_2|$ является первым интегралом системы, то $\dot{h}_2=0$. Из гамильтоновой системы~(\ref{eq:hamsys}) следует $h_3 u_1=0$. Случай $h_3=0$ противоречит условию нетривиальности начального ковектора в ПМП. В случае $u_1=0$ выражения~(\ref{eq:abnormh1min}) и ~(\ref{eq:abnormh1eq}) совпадают. Таким образом, экстремальное управление при $h_1=0$ имеет тот же вид, что и при $h_1<0$.

При $h_1>0$ из~(\ref{eq:HamMax}) следует $\sqrt{h_1^2+h_2^2}=0 \Leftrightarrow h_1=h_2=0$, что противоречит предположению $h_1>0$. Таким образом, этот случай невозможен. При $h_1>0$ анормальных экстремалей не существует.

Таким образом, мы получили следующее утверждение
\begin{theorem}\label{thm:abnormal}
Анормальное управление существует при $h_1\leq0$ и имеет вид $u_1(t)=0$, $u_2(t) \in I=[-1,1]$ --- произвольная функция класса $L_\infty([0,T],I)$, удовлетворяющая условию $h_{10} \cos U_2(t) - h_{30} \sin U_2(t)<0$ при всех $t \in [0,T]$. 
\end{theorem}

Анормальные траектории --- решения горизонтальной части системы~(\ref{eq:hamsysh1na}). Они являются поворотами с постоянной скоростью $s_2=\pm 1$:
\begin{equation*}\label{eq:exttrajdbnorm}
x(t)=0,\quad y(t)=0,\quad \theta(t)=s_2 t. 
\end{equation*}

\subsection{Нормальный случай $H=1$}
Вертикальная часть гамильтоновой системы ПМП~(\ref{eq:hamsys}) имеет вид 
\begin{equation}\label{eq:vertpartsys}
\dot{h}_1 = -u_2 h_3, \qquad
\dot{h}_2 = u_1 h_3, \qquad
\dot{h}_3 =u_2 h_1,
\end{equation}
где управление $(u_1, u_2)$ задано соотношениями~(\ref{eq:abnormh1min}) --(\ref{eq:abnormh1max}).

Система~(\ref{eq:vertpartsys}) имеет первый интеграл --- гамильтониан $H$, заданный формулой~(\ref{eq:HamMax}), зависящей от знака $h_1$. Помимо гамильтониана имеется еще один первый интеграл --- функция Казимира~$E$, заданная единственной формулой на всей области определения: 
\begin{equation}\label{eq:Casimir}
E = h_1^2 + h_3^2.
\end{equation}
\begin{remark}\label{rem:Casimir}
Функции Казимира --- это функции на двойственном пространстве алгебры Ли, коммутирующие в смысле скобок Пуассона со всеми левоинвариантными гамильтонианами. Они являются универсальными законами сохранения на группе Ли.
 Связные совместные поверхности уровня всех функций Казимира являются орбитами коприсоединенного представления группы Ли (см.~\cite[Утв.~7.7]{Laurent-Gengoux}). 

Первый интеграл $E$ является функцией Казимира на $\SEtwo$, поскольку он коммутирует со всеми левоинвариантными гамильтонианами:
$$\{E,h_1\}=\{E,h_2\}=\{E,h_3\}=0.$$ 
\end{remark}

\begin{figure}[ht]
\centering
\includegraphics[width=\hsize]{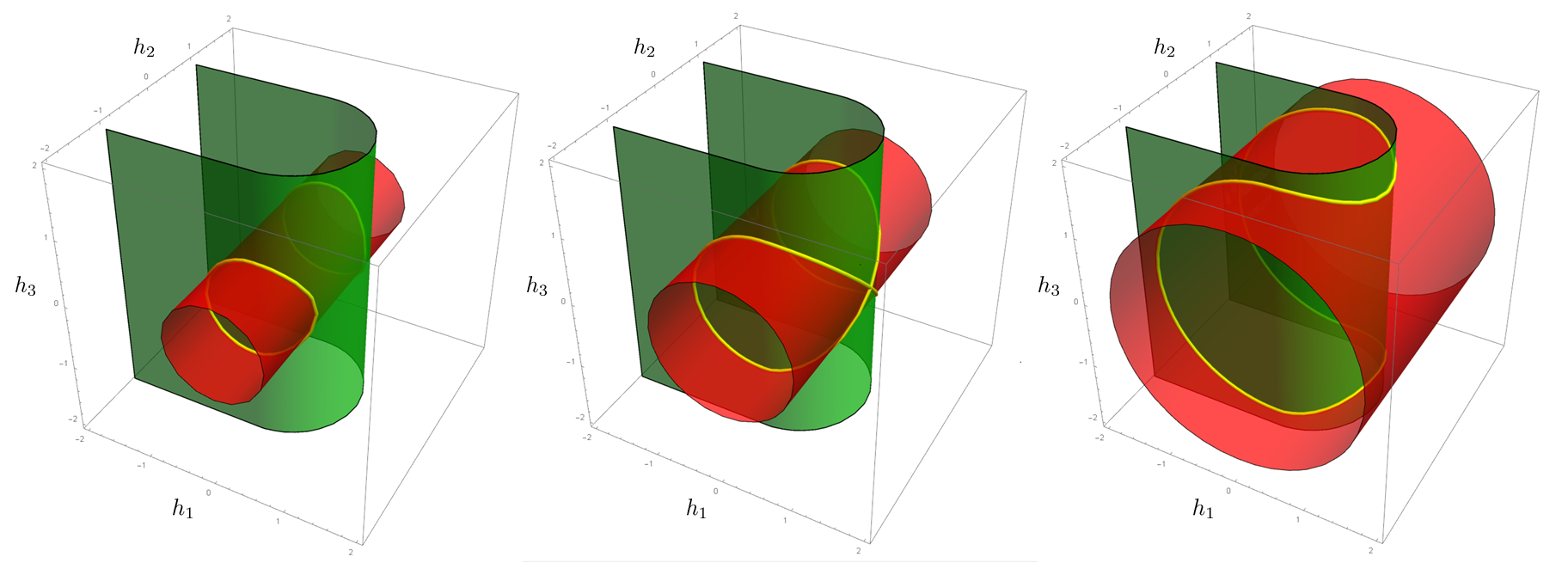}
\caption{Взаимное расположение поверхностей уровня гамильтониана $H$ (зеленая поверхность) и функции Казимира $E$ (красная поверхность). Линия пересечения выделена желтым. Слева: $E<1$. Центр: $E=1$. Справа: $E>1$.}
\label{fig:phaseplot3D}
\end{figure}

На Рис.~\ref{fig:phaseplot3D} изображены различные варианты взаимного расположения поверхности уровня гамильтонина $H=1$, представляющей собой две полуплоскости $\{(h_1,h_3)\;|\; h_1<0, \, h_2 = \pm 1\}$, склеенные половиной цилиндра $\{h_1^2+h_2^2=1\; | \; h_1\geq 0\}$, и поверхности уровня функции Казимира $E \geq 0$, представляющей собой цилиндр $h_1^2+h_3^2=E$. Помимо вырожденных случаев $\{h_1=0,\, h_2 =\pm 1,\, h_3=0\}$ (устойчивые положения равновесия) и $\{h_1=1,\, h_2 =0, \, h_3=0\}$ (неустойчивое положение равновесия), возможны три различных случая: $E<1$, $E=1$ и $E>1$. Случай $E=1$ называется критическим. В этом случае движение в состояние (или из состояния) неустойчивого равновесия происходит за бесконечное время. 

\begin{figure}[ht]
\centering
\includegraphics[width=0.9\hsize]{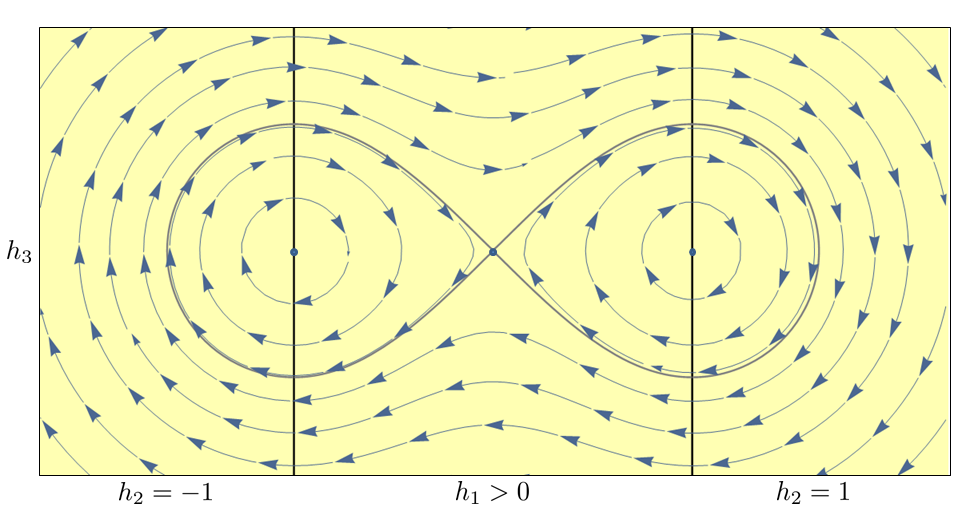}
\caption{Фазовый портрет на поверхности уровня гамильтониана $H=1$.}
\label{fig:phaseplotFold}
\end{figure}

На Рис.~\ref{fig:phaseplotFold} приведен фазовый портрет вертикальной части гамильтоновой системы, суженный на поверхность уровня гамильтониана $H=1$. На этом рисунке поверхность $H=1$, состоящая из двух полуплоскостей, склеенных половиной цилиндра, изображена в развернутом виде. Линии склейки нанесены на рисунок. Левая часть рисунка соответствует фазовому портрету на полуплоскости $h_2=-1$, средняя часть ---  на половине цилиндра, правая часть ---  на полуплоскости $h_2=1$. По вертикали обозначены значения $h_3$. По горизонтали --- значения $h_1$ на левой и правых частях рисунка, значение $\arccos h_1 \in [-\frac{\pi}{2},\frac{\pi}{2}]$ на средней части рисунка. Траектория, соответствующая критическому значению $E=1$ нанесена на рисунок.  

В зависимости от знака $h_1$ получаются две различные системы, см.~(\ref{eq:HamMax}). При смене знака $h_1$ динамика переключается с одной системы на другую. Далее мы проанализируем возможные случаи. Момент времени, в который происходит переключение, будем обозначать $t_0 \in \{t_0^0=0,t_0^1,t_0^2,\ldots\}$. 

\subsubsection{Случай $h_1<0$}\label{subsec:h1min}
Заметим, что в силу постоянства гамильтониана $H$, см.~(\ref{eq:HamMax}), выполнено $$h_2=h_{20}=\pm 1.$$

В силу~(\ref{eq:abnormh1min}) имеем $u_1=0$, $u_2 = h_2$.
Обозначим $s_2=h_2=\pm1$. 

Гамильтонова система ПМП имеет вид
\begin{equation}\label{eq:hamsysh1n}
\left\{\begin{array}{l l}
\dot{x} = 0,  & x(t_0)=x_0,\\[0.005\linewidth]
\dot{y} =0,  & y(t_0)=y_0,\\[0.005\linewidth]
\dot{\theta} =s_2, & \theta(t_0)=\theta_0,
\end{array}\right.
\quad
\left\{\begin{array}{l l}
\dot{h}_1 = -s_2 h_3,  & h_1(t_0)=h_{10},\\[0.005\linewidth]
\dot{h}_3 =s_2 h_1,& h_3(t_0)=h_{30}.
\end{array}\right.
\end{equation}

Решение горизонтальной части имеет вид
\begin{equation}\label{eq:exttrajh1n}
x(t)=x_0,\quad y(t)=y_0,\quad \theta(t)=\theta_0+s_2 (t-t_0). 
\end{equation}

Получаем, что экстремальными траекториями являются повороты плоскости вокруг неподвижной точки $x_0$, $y_0$ с постоянной скоростью $s_2=\pm 1$.


Вертикальная часть задана линейной системой ОДУ, которая легко интегрируется:
\begin{equation}\label{eq:h1h3}
\left\{\begin{array}{l l}
h_1(t) =  h_{10} \cos(t-t_0)-s_2 \,h_{30} \sin(t-t_0),\\[0.005\linewidth]
h_3(t) =  h_{30} \cos(t-t_0)+s_2 \,h_{10} \sin(t-t_0).
\end{array}\right.
\end{equation}

Получаем, что решением вертикальной части являются дуги окружностей в плоскости $(h_1, h_3)$, а если точнее, то в полуплоскости $h_1<0$ в силу начального предположения. При $s_2=-1$ движение по окружности происходит по часовой стрелке, а при $s_2=1$ --- против часовой стрелки.

Для полного исследования рассматриваемого случая осталось вычислить момент времени $t_1$, в который произойдет переключение динамики, то есть такой момент, когда условие $h_1<0$ перестанет выполняться
$$t_1=\min\left\{t>t_0\; | \; h_{10} \cos(t-t_0)-s_2 h_{30} \sin(t-t_0) =0\right\}.$$

Геометрически это уравнение означает, что $t_1-t_0$ --- минимальный угол, на который нужно повернуть плоскость $(h_1,h_3)$, чтобы вектор $(1,0)$ стал перпендикулярным вектору $v=(h_{10}, - s_2 h_{30})$. Рассматривая два возможных варианта расположения $v$ 
в полуплоскости $h_1<0$, находим
\begin{equation}\label{eq:tswitchn}
t_1-t_0=\arg\left(-s_2 h_{30} - \i h_{10}\right) \in (0,\pi].
\end{equation}

В заключение отметим, что в случае $t_0>0$, то есть когда уже произошло хотя бы одно переключение, формула~(\ref{eq:tswitchn}) сводится к $t_1 - t_0 = \pi$.  

\subsubsection{Случай $h_1=0$}\label{subsec:h1eq}

В силу постоянства гамильтониана $H$, см.~(\ref{eq:HamMax}), выполнено $h_2=h_{20}=\pm 1.$ 
В силу~(\ref{eq:abnormh1eq}) имеем $u_1=0$, $u_2 = h_2$.
Обозначим $s_2=h_2=\pm1$. \\
Динамика $h_1$ зависит от знака производной $\dot{h}_1= -h_2 h_3$. 

В случае $h_{30}=0$ система находится в устойчивом стационарном состоянии $h_1=h_3=0$, $h_2=s_2$. Соответствующее экстремальное управление постоянно $u_1=0$, $u_2 =s_2$, а экстремальной траекторией является вращение на месте $x(t)=x_0$, $y(t)=y_0$, $\theta(t)=\theta_0+s_2 (t-t_0)$. Очевидно, что такие траектории являются оптимальными при $t-t_0 \in [0, \pi]$.

В случае  $h_{30}\neq 0$ обозначим $s_{23}=\sign h_2 h_3$. Тогда $\dot{h}_1=-s_{23}=\pm 1$. Если $s_{23}=1$, то динамика задается системой~(\ref{eq:hamsysh1n}) рассмотренного ранее случая $h_{1}<0$. Если $s_{23}=-1$, то динамика задается системой~(\ref{eq:hamsysh1p})  рассматриваемого далее случая $h_{1}>0$. 

\subsubsection{Случай $h_1>0$}\label{subsec:h1max}
В силу~(\ref{eq:HamMax}) имеем $h_1=\sqrt{1-h_2^2}$. Гамильтонова система ПМП имеет вид
\begin{equation}\label{eq:hamsysh1p}
\left\{\begin{array}{l l}
\dot{x} = h_1 \cos\theta,  & x(t_0)=x_0,\\[0.005\linewidth]
\dot{y} =h_1 \sin\theta,  & y(t_0)=y_0,\\[0.005\linewidth]
\dot{\theta} =h_2, & \theta(t_0)=\theta_0,
\end{array}\right.
\quad
\left\{\begin{array}{l l}
\dot{h}_1 = -h_2 h_3,  & h_1(t_0)=h_{10},\\[0.005\linewidth]
\dot{h}_2 = h_1 h_3, & h_2(t_0)=h_{20},\\[0.005\linewidth]
\dot{h}_3 =h_2 h_1,& h_3(t_0)=h_{30}.
\end{array}\right.
\end{equation}
Эта система является модельным примером в геометрической теории управления~\cite{notes}. Явное выражение решения в эллиптических функциях Якоби получено в статье~\cite{yuriSE2}, где авторы свели вертикальную подсистему к уравнению математического маятника и проинтегрировали его в терминах выпрямляющих координат. При этом решение задается различными 
 формулами в разных областях фазового портрета маятника. Конкретный вид формулы определяется характером движения маятника: колебание, вращение, движение по сепаратрисе, устойчивое или неустойчивое равновесие.

В данном разделе мы предлагаем другую технику интегрирования, ведущую к явной параметризации решения одной формулой почти всюду. Техника состоит в следующем: сначала мы выведем ОДУ на функцию $h_2$, 
найдем его явное решение, и затем выразим оставшиеся компоненты $h_1$, $h_3$ через известную функцию $h_2$ и начальные условия.  

Напомним, что помимо первого интеграла $h_1^2+h_2^2=1$, имеется еще один первый интеграл --- функция Казимира $E=h_1^2+h_3^2$. Обозначим $M=E-2$.

Выписывая вторую производную функции $h_2$ в силу~(\ref{eq:hamsysh1p}), получаем ОДУ
\begin{equation}\label{eq:ODEh2}
\ddot{h}_2+M h_2+2 h_2^3=0,
\end{equation} 
с начальными условиями
\begin{equation}\label{eq:h2init}
h_2(t_0) = h_{20},\quad \dot{h}_2(t_0)=h_{10}\, h_{30}=\sqrt{1-h_{20}^2} \, h_{30}.
\end{equation} 

\paragraph{Случай $h_{30} =0$} В силу~(\ref{eq:HamMax}) имеем $|h_{20}|<1$. При  $h_{20}=0$ решением системы~(\ref{eq:ODEh2})--(\ref{eq:h2init}) является неустойчивое положение равновесия. Соответствующей экстремальной траекторией $\gamma(t)=(x_0+(t-t_0) \cos\theta_0,y_0+(t-t_0) \sin\theta_0, \theta_0)$ является движение по прямой. Такие траектории являются оптимальными до бесконечности. При $h_{20}\neq 0$ система сводится к рассматриваемому далее случаю $h_{30}\neq 0$ с углом $\alpha = \frac{\pi}{2} \sign h_{20} $.

\paragraph{Случай $h_{30} \neq 0$} Рассмотрим далее общий случай $E\not\in\{0,1\}$. Частные случаи $E\in\{0,1\}$ выражаются в элементарных функциях и получаются предельным переходом в формулах при $E<1$ и $E>1$. 

Обозначим $s_2=\sign h_{20}$, $s_3 = \sign h_{30}$. Введем угол
$$\alpha = \left\{\begin{array}{l l}
\arg\left(-s_3 \left(h_{20} + \i \sqrt{1-h_{20}^2}\right)\right) \in (-\pi, \pi], & \text{ при } E>1,\\[0.02\linewidth]
\arg \left(h_{20} + \i \sqrt{1-h_{20}^2}\right)+\frac{1-s_2}{2} \pi \in (-\frac{\pi}{2}, \frac{\pi}{2}], & \text{ при } E<1.
\end{array} \right. $$
Обозначим
$$\xi(t) = \left\{\begin{array}{l l}
\frac{t-t_0}{k}-s_3 F(\alpha,k), & \text{ при } E>1,\\[0.02\linewidth]
-s_2 s_3 \frac{t-t_0}{k}+s_2 F(\alpha,k),  & \text{ при } E<1,
\end{array} \right.
$$
где $F$ --- эллиптический интеграл 1-го рода в форме Лежандра:
$$F(\alpha,k)=\int\limits_0^\alpha \frac{\dd a}{\sqrt{1-k^2 \sin^2 a}}.$$

Явное решение уравнения~(\ref{eq:ODEh2}) для общих граничных условий приведено в~\cite[App.~A]{Lakshmanan}, см. также~\cite{Mathews}. Подставляя граничные условия~(\ref{eq:h2init}), получаем
\begin{equation}\label{eq:h2sol}
h_2(t) = -s \cn\left(\xi(t),k\right), 
\end{equation} 
где $$k=\frac{1}{\sqrt{E}}=\frac{1}{\sqrt{1-h_{20}^2+h_{30}^2}}, \qquad 
s=\left\{\begin{array}{l l} 
s_3, & \text{ при } E>1,\\
-s_2, & \text{ при } E<1.
\end{array}\right.$$


Из~(\ref{eq:h2sol}) и условия $h_1(t)=\sqrt{1-h_2^2(t)}$ находим
\begin{equation}\label{eq:solh1}
h_1(t) = \sn\left(\xi(t),k\right), 
\end{equation}

Функция $h_3(t)$ находится из системы $\dot{h}_3(t)=h_1(t)h_2(t)$, $h_3(0)=h_{30}$:
\begin{equation*}\label{eq:solh3deriv}
\begin{array}{r l}
h_3(t) = & \displaystyle h_{30}+\int\limits_{t_0}^t h_1(\tau)h_2(\tau) \dd \tau = h_{30}+\frac{s}{k} \int\limits_{\xi(t_0)}^{\xi(t)}-k^2 \sn\left(a,k\right) \cn\left(a,k\right) \dd a.
\end{array}
\end{equation*}

Используя тождество $\frac{d}{d a}\dn(a,k)=-k^2 \sn(a,k) \cn(a,k)$, получаем
\begin{equation}\label{eq:solh3}
h_3(t) = h_{30}+\frac{s}{k}\left(\dn\left(\xi(t),k\right) - \dn\left(\xi(t_0),k\right)\right).
\end{equation}


Обозначим
\begin{equation}\label{eq:H2}
\displaystyle H_2(t)=\int\limits_{t_0}^{t} h_2(\tau) \dd\, \tau.
\end{equation}

\begin{remark}\label{rem:h1h3}
По аналогии с результатом~\cite[Гл.~1,~$\S$~3]{arnold} задача Коши
\begin{equation*}
\label{eq:Cauchysyst}
\begin{cases}
\dot{h}_1 = -h_2 h_3,  & h_1(t_0)=h_{10},\\[0.005\linewidth]
\dot{h}_3 =h_2 h_1,& h_3(t_0)=h_{30}.
\end{cases}
\end{equation*}
имеет единственное решение $(h_1,h_3)$,  заданное формулой
\begin{equation*}\label{h1h3sol}
\begin{array}{l}
h_1(t)=h_{10}\cos H_2(t)-h_{30} \sin H_2(t),\\[0.005\linewidth]
h_3(t)=h_{30}\cos H_2(t)+h_{10} \sin H_2(t).
\end{array}
\end{equation*}
Заключаем, что как и при $h_1<0$ проекции решений вертикальной части на плоскость $(h_1,h_3)$ являются дугами окружностей, однако в отличие от случая $h_1<0$, движение по ним происходит с переменной скоростью.
\end{remark}

Интеграл~(\ref{eq:H2}), где $h_2$ задается формулой~(\ref{eq:h2sol}), допускает явное представление в эллиптических функциях Якоби, см.~\cite{Handbook},
\begin{equation}\label{eq:H2exact}
H_2(t)=-s\arccos\left(\dn\left(\xi(\tau),k\right)\right)\Big|_{\tau=t_0}^{\tau=t}.
\end{equation}

Теперь запишем решение горизонтальной части гамильтоновой системы. Заметим, что в силу $\dot{\theta}=h_2$, функции $\theta$ и $H_2$, см.~(\ref{eq:H2}), связаны соотношением
\begin{equation}\label{eq:thh1p}
\theta(t)=\theta_0+H_2(t).
\end{equation}

Компоненты $x$, $y$ экстремальной траектории выражаются в квадратурах
\begin{equation}\label{xysol}
x(t)=x_0+\int\limits_{t_0}^t h_1(\tau) \cos\theta(\tau) \, \dd \tau,
\quad
y(t)=y_0+\int\limits_{t_0}^t h_1(\tau) \sin\theta(\tau) \, \dd \tau.
\end{equation}

Введем систему координат $(\tilde{x},\tilde{y},\tilde{\theta})$ в окрестности точки $(x_0,y_0,\theta_0)$:
$$
\left(\begin{array}{c}
\tx \\
\ty
\end{array}\right)=\left(\begin{array}{c c}
\cos \beta_0 & \sin\beta_0 \\
-\sin \beta_0 & \cos\beta_0
\end{array}\right)
\left(\begin{array}{c}
x \\
y
\end{array}\right) -
\left(\begin{array}{c}
x_0 \\
y_0
\end{array}\right),
\quad
\tth = \theta - \beta_0.
$$
где $$\beta_0=\theta_0+s \arccos\left(\dn\left(\xi(t_0),k\right)\right).$$

В силу левоинвариантности гамильтоновой системы, компоненты траектории $(\tx(t),\ty(t),\tth(t))$ в новых координатах удовлетворяют тем же уравнениям~(\ref{eq:hamsysh1p}), что и в исходных координатах:
$$\dot{\tx}(t) = h_1(t) \cos \tth(t),\quad \dot{\ty}(t) = h_1(t) \sin \tth(t), \quad \dot{\tth}(t)=h_2(t).$$

Из тождества~(\ref{eq:thh1p}) следует
\begin{equation}\label{eq:tth}
\tth(t) = -s\arccos\left(\dn\left(\xi(t),k\right)\right).
\end{equation}

Компоненты $\tx$, $\ty$ экстремальной траектории выражаются как
\begin{equation}\label{txysol}
\tx(t)=\int\limits_{t_0}^t h_1(\tau) \cos\tth(\tau) \, \dd \tau,
\quad
\ty(t)=\int\limits_{t_0}^t h_1(\tau) \sin\tth(\tau) \, \dd \tau.
\end{equation}

Интегралы в правой части допускают явное представление. Действительно, подставляя~(\ref{eq:solh1}), (\ref{eq:tth}) в (\ref{txysol}), имеем
\begin{equation*}
\begin{array}{r l}
\tx(t)=& \displaystyle \int\limits_{t_0}^t \sn\left(\xi(\tau),k\right) \dn\left(\xi(\tau),k\right) \, \dd \tau = k \int\limits_{\xi(t_0)}^{\xi(t)}\sn\left(a,k\right) \dn\left(a,k\right) \, \dd a\\
=& -k \left(\cn\left(\xi(t),k\right)-\cn\left(\xi(t_0),k\right)\right),
\\
\ty(t)=&  \displaystyle-s \int\limits_{t_0}^t \sn\left(\xi(\tau),k\right) \sqrt{1-\dn^2\left(\xi(\tau),k\right)} \, \dd \tau =-s \, k^2 \int\limits_{\xi(t_0)}^{\xi(t)} \sn^2\left(a,k\right) \, \dd a\\
=& \displaystyle -s \left(a-\EE\left(a,k\right)\right)\Bigg|_{\xi(t_0)}^{\xi(t)}=
-s \left(\frac{t-t_0}{k}-\EE\left(\xi(t),k\right)+\EE\left(\xi(t_0),k\right)\right),
\end{array}
\end{equation*}
где $\EE(a,k)=E(\am(a,k),k)$ --- эллиптический интеграл второго рода:
$$E(\alpha,k) = \int\limits_0^\alpha \sqrt{1-k^2 \sin^2 a} \dd a.$$

При интегрировании выражения для $\ty(t)$ использовалось стандартное тождество $\sqrt{1-\dn^2(a,k)}=k \,|\sn(a, k)|$ и условие $\sn(\xi(t), k)>0$, выполненное в силу предположения $h_1(t)>0$, см.~(\ref{eq:solh1}). 

Возвращаясь к исходным координатам $(x,y,\theta)$, получаем
\begin{equation}\label{eq:xythfin}
\left(\begin{array}{c}
x \\
y
\end{array}\right)=\left(\begin{array}{c c}
\cos \beta_0 & -\sin\beta_0 \\
\sin \beta_0 & \cos\beta_0
\end{array}\right)
\left(\begin{array}{c}
\tx \\
\ty
\end{array}\right) +
\left(\begin{array}{c}
x_0 \\
y_0
\end{array}\right),
\quad
\theta=\tth + \beta_0.
\end{equation}

\begin{remark}\label{rem:sparam}
Заметим, что поскольку $h_1>0$, то траектория может быть параметризована длиной $s$ ее проекции на плоскость $(x,y)$: 
$$\s(t)=\int_{t_0}^t \sqrt{\dot{x}^2(\tau)+\dot{y}^2(\tau)} \dd \tau =\int_{t_0}^t h_1(\tau) \dd \tau.$$
Этот интеграл выражается явно, см.~\cite{Handbook}. В силу~(\ref{eq:solh1}) имеем
$$
\s(t) = \int\limits_{t_0}^t \sn\left(\xi(\tau),k\right) \dd \tau
 =  \ln\left(\dn\left(\xi(\tau),k\right) -k \cn\left(\xi(\tau),k\right)\right)\Big|_{\tau=t_0}^{\tau=t}.
$$
Явные выражения для компонент траектории $x(\s)$, $y(\s)$, $\theta(\s)$, параметризованной длиной дуги $\s$ на плоскости, найдено в~\cite{DuitsJMIV2014}.
\end{remark}

\subsubsection{Нормальные экстремальные траектории и управления}\label{seq:Extremals}
Подведем итог этого раздела, сформулировав следующую теорему
\begin{theorem}\label{thm:extremal}
В задаче быстродействия~(\ref{eq:contsyst})--(\ref{eq:bounds}) нормальное экстремальное управление $(u_1(t), u_2(t))$ однозначно определяется значением
$$h_0^0 = (h_{10}^0,s_{20}^0,h_{30}^0), \quad h_{10}^0 \in  (-\infty,1], \; s_{20}^0 \in \{-1,1\}, \; h_{30}^0 \in \R.$$

Функция~$u_1(t)$ имеет вид $u_1(t)= \sqrt{1-u_2^2(t)}$, $t \in [0,T]$. 

Функция~$u_2(t)=h_2(t)$ задается на интервалах, образованных разбиением отрезка $t\in [0,T]$ точками переключения $t_0 \in \{0=t_0^0, t_0^1,t_0^2, \ldots,T\}$, где момент переключения $t_0^i$ зависит от состояния $h_0^{i-1}=(h_{10}^{i-1},s_{20}^{i-1},h_{30}^{i-1})$, достигнутого
в предыдущий момент $t_0^{i-1}$, и
определяется рекуррентной формулой
$$t_0^{i}(h_0^{i-1}) = \min\{t>t_0^{i-1} \, | \, h_1(t,h_0^{i-1})=0\},\quad h_0^i=h(t_0^i(h_0^{i-1}), h_0^{i-1}).$$
Здесь $h(t,h_0^{i}) =(h_1(t,h_0^{i}),\, h_2(t,h_0^{i}),\, h_3(t,h_0^{i}))= \operatorname{vert}\left(e^{t \vec{H}}  h_0^{i} \right)$ --- решение вертикальной части гамильтоновой системы ПМП с начальным значением $h_0^i$ за время $t \geq t_0^i$, имеющей вид
$$
\left\{\begin{array}{l l} 
(\ref{eq:hamsysh1n}), & \text{ при } (h_{10}^i<0) \lor (h_{10}^i=0 \land s_{20}^i s_{30}^i>0),\\ 
(\ref{eq:hamsysh1p}), & \text{ при } (h_{10}^i>0) \lor (h_{10}^i=0 \land s_{20}^i s_{30}^i<0),
\end{array}\right. 
$$
где $s_{30}^i=\left\{\begin{array}{l l}
\sign h_{30}^i, & \text{ при } h_{30}^i\neq0,\\ 
s_{20}^i, & \text{ при }  h_{30}=0.\end{array}\right.$ 
\end{theorem}

Экстремальные траектории имеют вид
\begin{equation*}\label{eq:exttaj}
x(t)=\int\limits_{0}^t u_1(\tau) \cos\theta(\tau) \, \dd \tau,\;
y(t)=\int\limits_{0}^t u_1(\tau) \sin\theta(\tau) \, \dd \tau,\;
\theta(t)=\int\limits_{0}^t u_2(\tau) \, \dd \tau.
\end{equation*}

Экстремальное управление и траектория являются решением $e^{t \vec{H}}  h_0^{0}$ вертикальной и горизонтальной части гамильтоновой системы ПМП $\dot{h}=\vec{H}$ с начальным условием $h(0)=h_0^0$.
Явные формулы в общем случае $E \not\in \{0,1\}$ при $h_1<0$, $h_1=0$ и $h_1>0$ приведены в разделах~\ref{subsec:h1min},~\ref{subsec:h1eq} и ~\ref{subsec:h1max} соответственно. Случаю $h_1<0$ соответствуют отрезки экстремальных траекторий, являющиеся вращениями на месте. Случаю $h_1>0$ соответствуют отрезки экстремальных траекторий, являющиеся субримановыми геодезическими на $\SEtwo$, проекция на плоскость которых не имеет точек возврата.  Случай $h_1=0$ примыкает к случаю $h_1<0$, при $h_{20}^i h_{30}^i \geq 0$, и к случаю $h_1>0$, при $h_{20}^i h_{30}^i< 0$.

Решения в частных случаях $E \in \{0,1\}$ выражаются в элементарных функциях и получаются предельным переходом в формулах в общем случае. Они соответствуют одному из типов движения:
\begin{itemize}
\item $E=0 \Rightarrow h_{10} = 0$ --- поворот на месте;
\item $E = 1$, $h_{30}=0$ --- движение по прямой;
\item $E=1$, $h_{30}\neq 0$, $h_{10}>0$ --- движение по трактрисе,~см.~\cite{yuriSE2};
\item $E=1$, $h_{30}\neq 0$, $h_{10}<0$ --- поворот на месте;
\item $E=1$, $h_{10}=0$ --- поворот на месте или движение по трактрисе.
\end{itemize}

\begin{figure}[ht]
\centering
\includegraphics[width=0.49\linewidth]{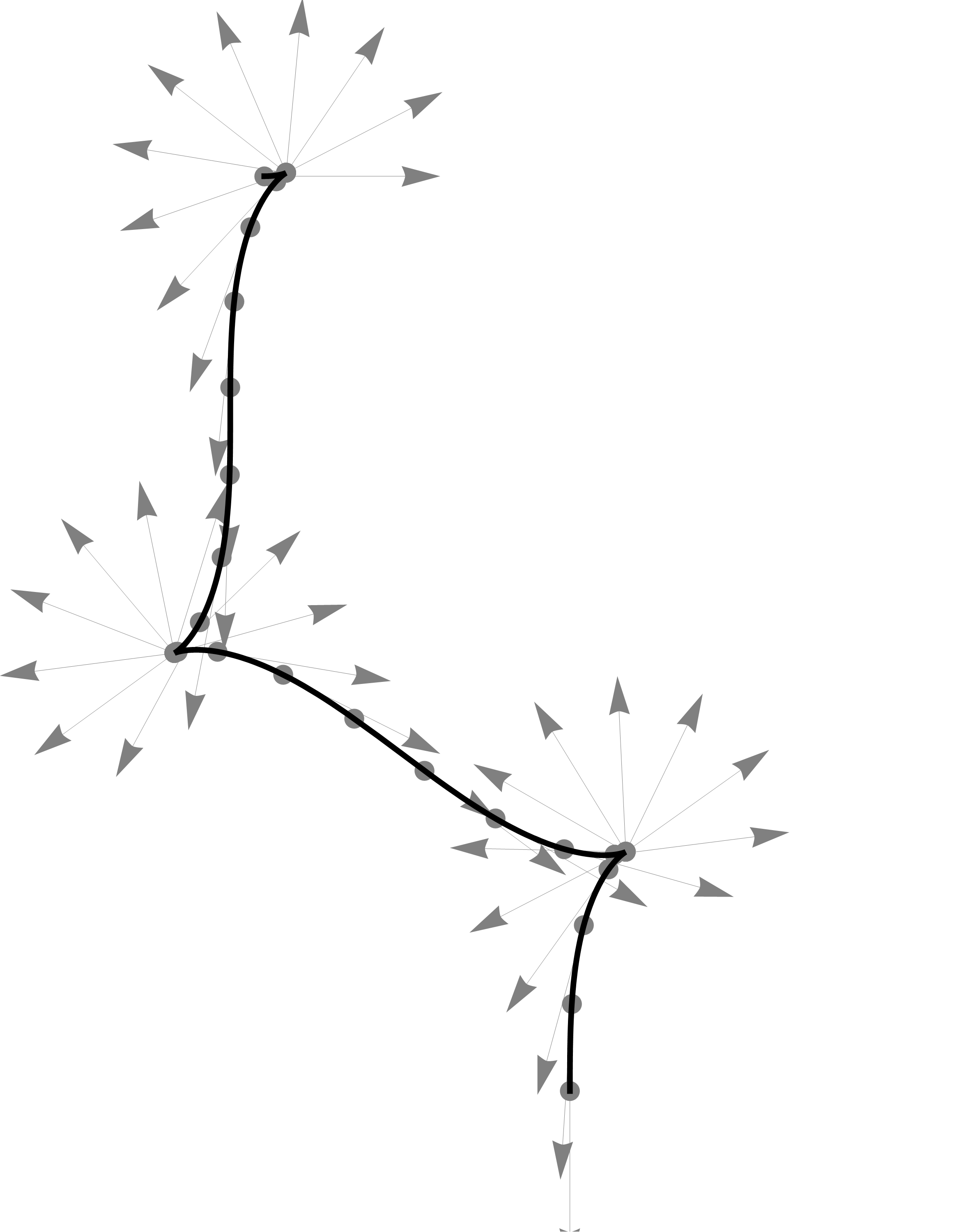}~
\includegraphics[width=0.49\linewidth]{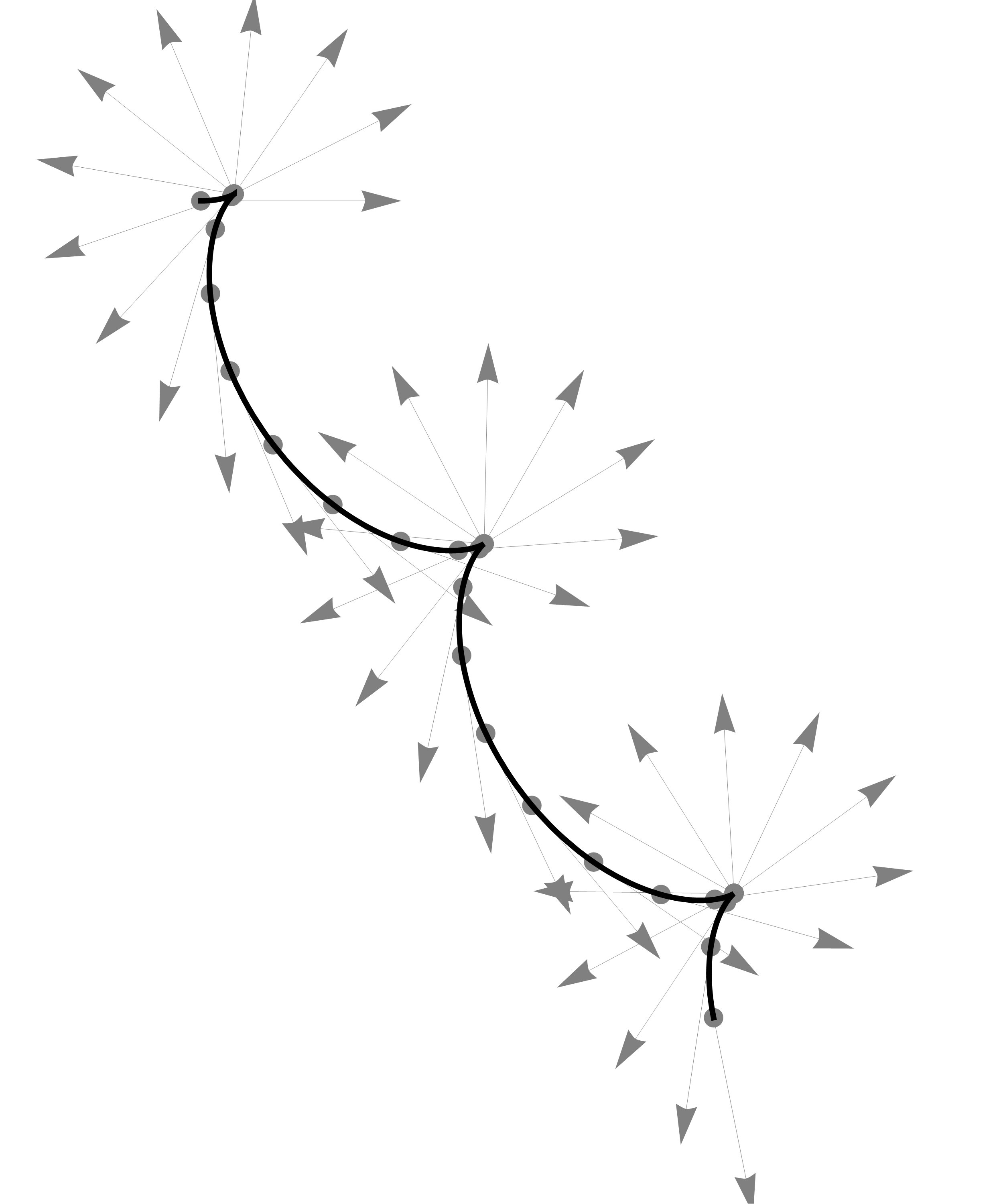}
\caption{Проекция экстремальных траекторий на плоскость $(x,y)$. Серой стрелкой обозначено направление $\theta$ в моменты времени $\{0,0.5,1,\ldots,20\}$. \textbf{Слева: } $h_{10}=\frac12$, $h_{20}=\frac{\sqrt{3}}{2}$, $h_{30}=1$. \textbf{Справа: }$h_{10}=\frac12$, $h_{20}=\frac{\sqrt{3}}{2}$, $h_{30}=0.7$.}
\label{fig:trajexamps}.
\end{figure}

Анализируя полученное решение, наблюдаем связь между экстремальными траекториями в исследуемой задаче  и задаче быстродействия для управляемой системы на группе движений плоскости с управлением в круге, исследованной в~\cite{yuriSE2}: проекции на плоскость $(x,y)$ траекторий в задаче на полукруге и задаче на круге совпадают, однако динамика угла $\theta$ при движении по ним отличается --- в точках возврата в задаче на полукруге происходит равномерное увеличение угла $\theta$ на $\pi$ радиан, см. Рис.~\ref{fig:trajexamps}.




\section{Оптимальный синтез}\label{sec:optimal}
В этом разделе произведем анализ оптимальности экстремальных управлений и опишем структуру оптимального синтеза.

Сначала заметим, что в задаче~(\ref{eq:contsyst})--(\ref{eq:bounds}) нет оптимальных строго анормальных траекторий. То есть любая оптимальная анормальная траектория является также нормальной. Это утверждение следует из Теоремы~\ref{thm:abnormal} и очевидного соображения, что поворот на месте оптимален (происходит наискорейшим способом) тогда и только тогда, когда направление вращения не меняется, и скорость максимальна $|u_2|=1$. Таким образом, достаточно рассмотреть только нормальный случай.

В работе~\cite{Duits18} указано, что оптимальная траектория на плоскости не имеет внутренних точек разворота. Точками разворота могут быть лишь граничные точки. Это утверждение формализуется в виде следующей теоремы. 

\begin{theorem}\label{thm:optimality}
Оптимальная траектория в задаче~(\ref{eq:contsyst})--(\ref{eq:bounds}) имеет вид
\begin{equation}\label{eq:opttraj}
\textrm{
\begin{tabular}{ c | c c c }
$t\in         $&$ [0,t_0^1) $&$ [t_0^1, t_0^2)$&$ [t_0^2,T]$ \\
\hline
$x(t)         $&$ 0             $&$ x_s(t)              $&$ x_1$ \\
$y(t)         $&$ 0             $&$ y_s(t)              $&$y_1$ \\
$\theta(t)  $&$ s_1 t        $&$ \theta_s(t)      $&$\theta_s(t_0^2) + s_2 (t-t_0^2)$,
\end{tabular}
}
\end{equation}
где $0\leq t_0^1 \leq t_0^2 \leq T$ --- моменты переключения управления, знаки $s_i = \pm 1$ --- параметры, определяемые граничными условиями, а траектория 
$$(x_s(t),y_s(t),\theta_s(t))=:q_s(t), \quad q_s(t_0^1)=(0,0,\theta_0^1),\; q_s(t_0^2)=(x_1,y_1,\theta_0^2)$$
 является субримановой кратчайшей на $\SEtwo$, проекция на плоскость которой не содержит внутренних точек возврата (т.е. для любого $t\in (t_0^1, t_0^2)$ выполнено $\dot{x}_s(t)^2+ \dot{y}_s(t)^2>0$).
\end{theorem}
\textit{Доказательство.} 
Вид~(\ref{eq:opttraj}) оптимальной траектории следует из Теоремы~\ref{thm:extremal} с учетом того, что повороты на месте могут быть только на начальном или конечном интервале времени. 
Неоптимальность траектории $\gamma$ с внутренней точкой разворота можно доказать, построив более быструю траекторию $\gamma_0$ (срезку). См. Рис.~\ref{fig:optimal}.

Рассмотрим траекторию $\gamma$ с внутренней точкой разворота. Обозначим  через $\bar{\gamma}$ проекцию $\gamma$ на плоскость. Обозначим точку разворота через $B = \bar{\gamma}(t)$, $t \in [t_B,t_B+\pi]$. Для любых двух точек $A = \bar{\gamma}(t_A)$ и $C= \bar{\gamma}(t_C)$, $t_A<t_B<t_B+\pi<t_C$ из достаточно малой окрестности точки $B$ выполнено 
\begin{equation}\label{eq:thineq}
|\theta_0| + |\theta_1|<\pi.
\end{equation}

Время движения $T=t_C-t_A$ по $\bar{\gamma}$ из $A$ в $C$ имеет оценку снизу 
\begin{equation}\label{eq:lgam}
l_{AB}+l_{BC} + \pi < T, 
\end{equation}
где $l_{AB}$ --- евклидово расстояние между точками $A$ и $B$, $l_{BC}$ --- между точками $B$ и $C$. Эта оценка верна, поскольку время движения вдоль любой траектории между двумя точками не может быть меньше, чем при движении по отрезку прямой. Поэтому $l_{AB}+l_{BC}< t_B-t_A + t_C-t_B-\pi=T-\pi$

Для расстояния $l_{AC}$ между $A$ и $C$ выполнено неравенство треугольника 
\begin{equation}\label{eq:triagineq}
l_{AC} < l_{AB} + l_{BC}. 
\end{equation}
Траектория $\gamma_0$ строится из трех последовательных движений: поворот на месте на угол $\theta_0$, движение вперед за время $l_{AC}$, поворот на месте на угол~$\theta_1$. Время движения $T_0$ по траектории $\gamma_0$ из $\gamma(t_A)$ в $\gamma(t_C)$ вычисляется явно 
\begin{equation}\label{eq:T0}
T_0=|\theta_0| +l_{AC}+ |\theta_1|.
\end{equation}

Комбинируя~(\ref{eq:thineq})--(\ref{eq:T0}), получаем $T_0<T$. 
 $\blacksquare$

\begin{figure}[ht]\label{fig:optimal}
\centering
\includegraphics[scale=1]{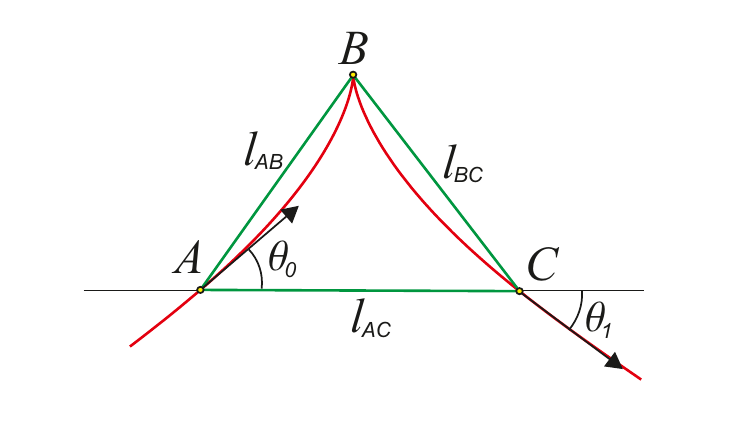}
\caption{Неоптимальность траектории с внутренней точкой поворота.}
\label{fig:optimal}
\end{figure}

Ввиду сложной параметризации экстремальных траекторий, в общем случае явное выражение оптимальной траектории $q_s$ через граничные условия в настоящее время остается неизвестным. Для этого требуется найти обратное экспоненциальное отображение (вообще говоря, многозначное) 
$$\Exp^{-1}: \M \to C \times \R^+: q_1 \mapsto (h(0), T),$$
 которое переводит конечную точку $q_1 \in \M$ оптимальной траектории в начальный импульс $h(0) \in C=T^*_e \M \cap \{H=1\}$ и время $T \geq 0$. В работе~\cite{DuitsJMIV2014} предложено численное решение, основанное на методе стрельбы.

\section*{Заключение}

В данной работе исследована задача быстродействия~(\ref{eq:contsyst})--(\ref{eq:bounds}) на группе движений плоскости с управлением в полукруге. Рассматриваемая управляемая система задает модель машины на плоскости, которая может двигаться вперед и поворачивать сколь угодно быстро. 

Траектории такой системы используются в области обработки изображений. Данная постановка нацелена на устранение проблемы точек возврата у существующего метода поиска выделяющихся кривых на изображениях. 

Помимо прикладного значения задача представляет интерес в теории управления, как модельный пример трехмерной управляемой системы с компактным множеством управлений, содержащим ноль на границе. 


В статье получены следующие основные результаты:
\begin{itemize}
\item конструктивно доказана полная управляемость и существование оптимальных управлений (Теорема~\ref{thm:exist});
\item найден явный вид анормальных управлений  (Теорема~\ref{thm:abnormal});
\item найдена явная параметризация нормальных управлений  (Теорема~\ref{thm:extremal});
\item найден явный вид экстремальных траекторий;
\item частично исследован вопрос оптимальности экстремалей;
\item описана структура оптимального синтеза (Теорема~\ref{thm:optimality}). 
\end{itemize}

Автор благодарен проф. Ю.\,Л. Сачкову за ценные замечания по работе.

 
\end{fulltext}

\end{document}